\documentclass[a4paper,10pt]{article}

\usepackage{amsfonts}
\usepackage{amsmath}
\usepackage{amssymb}

\usepackage{graphpap}
\usepackage{epsfig}

\title{Convergence of symmetric trap models in the hypercube}
\author{L.~R.~G.~Fontes\thanks{Partially
supported by CNPq grants 475833/2003-1,
307978/2004-4 and 484351/2006-0, and FAPESP grant 2004/07276-2.}\,\,\footnotemark[3]
\and P.~H.~S.~Lima\thanks{Supported by FAPESP grant 2004/13009-7}\,\,\thanks{Instituto de Matem\'atica
e Estat\'{\i}stica, Universidade de S\~ao Paulo, Rua do Mat\~ao
1010, Ci\-da\-de Universit\'aria, 05508-090 S\~ao Paulo  SP,
Brasil, emails: {\{lrenato, plima\}@ime.usp.br}}}
\date{}
\begin{document}
\maketitle

\def\Z{{\mathbb Z}}
\def\V{{\mathbb V}}
\def\N{{\mathbb N}^*}
\def\bN{\bar{\mathbb N}^*}
\def\H{{\mathcal H}}
\def\V{{\mathcal V}}
\def\ll{{\mathcal L}}
\def\K{{\mathcal K}}
\def\X{{\mathcal X}}
\def\zz{{\cal Z}}
\def\D{{\cal D}}
\def\M{{\cal M}}
\def\W{{\cal W}}
\def\J{{\cal J}}
\def\T{{\cal T}}
\def\O{{\cal O}}
\def\vt{\vec{\T}}
\def\vn{\vec{\N}}

\def\a{\alpha}
\def\b{\beta}
\def\ttau{\tilde\tau}
\def\htau{\hat\tau}

\def\g{\gamma}
\def\G{\Gamma}
\def\e{\epsilon}
\def\d{\delta}
\def\l{\lambda}
\def\tl{\tilde\lambda}
\def\L{\Lambda}
\def\tg{\tilde\gamma}
\def\tgd{\tilde\gamma^{d}}
\def\xd{X^{d}}
\def\hy{\hat Y^{d}}

\def\ve{\vec\eta}
\def\vz{{\vec\zeta}}

\def\P{{\mathbb P}}
\def\E{{\mathbb E}}
\def\pp{\P_p}
\def\ep{\E_p}
\def\C{{\mathcal C}}
\def\vc{\vec\C}

\def\s{\sigma}

\def\={&=&}

\newcommand{\stack}[2]{\genfrac{}{}{0pt}{3}{#1}{#2}}

\newtheorem{theo}{Theorem}
\newtheorem{prop}[theo]{Proposition}
\newtheorem{lm}[theo]{Lemma}
\newtheorem{cor}[theo]{Corollary}
\newtheorem{rmk}[theo]{Remark}
\newtheorem{df}[theo]{Definition}
\def\beq{\begin{equation}}
\def\eeq{\end{equation}}
\newcommand{\bei}{\begin{itemize}}
\newcommand{\eei}{\end{itemize}}
\newcommand{\ben}{\begin{enumerate}}
\newcommand{\een}{\end{enumerate}}
\newcommand{\beqn}{\begin{eqnarray}}
\newcommand{\beqnn}{\begin{eqnarray*}}
\newcommand{\eeqn}{\end{eqnarray}}
\newcommand{\eeqnn}{\end{eqnarray*}}
\newcommand{\brm}{\begin{rmk}}
\newcommand{\erm}{\end{rmk}}
\renewcommand{\theequation}{\thesection .\arabic{equation}}

\def\xd{X^{d}}


\begin{abstract}
We consider symmetric trap models in the $d$-dimensional hypercube whose ordered mean waiting times,
seen as weights of a measure in $\N$, converge to a finite measure as $d\to\infty$, and show that
the models suitably represented converge to a K process as $d\to\infty$. We then apply this 
result to get K processes as the scaling limits of the REM-like trap model and the Random Hopping 
Times dynamics for the Random Energy Model in the hypercube in time scales corresponding to the 
{\em ergodic} regime for these dynamics.
\end{abstract}

\paragraph{Keywords and phrases} trap models, weak convergence, K-process, random energy model, scaling limit, aging

\paragraph{2000 Mathematics Subject Classification} 60K35, 60K37, 82C44

\section{Introduction}
\setcounter{equation}{0}
\label{sec:intro}

Trap models have been proposed as qualitative models exhibiting localization and aging
(see~\cite{kn:NE, kn:B} for early references).
In the mathematics literature there has recently been an interest in establishing
such results for a varied class of such models (see~\cite{kn:FIN,kn:ABG,kn:ACM} and references therein).
In particular, it has been recognized that scaling limits play an important role in such
derivations (see~\cite{kn:FIN,kn:AC1,kn:FM,kn:ABC} and references therein). It may be argued
that such phenomena correspond to related phenomena exhibited by limiting models.

In this paper we consider symmetric trap models in the hypercube whose mean waiting times
converge as a measure to a finite measure as the dimension diverges, and show that these
models converge weakly. We then apply this result to establish the scaling limits of two
dynamics in the hypercube, namely the REM-like trap model and the Random Hopping Times dynamics
for the Random Energy Model, in time scales corresponding to the ergodic regime for these dynamics.

\subsection{The model}
\setcounter{equation}{0}
\label{ssec:modres}

Let $\mathcal{H}$ denote the $d$-dimensional hypercube, namely $\mathcal{H}$ is the graph
$(\mathcal{V},\mathcal{E})$ with
\[
\begin{array}{l}
\mathcal{V} =  \{0,1\}^d,\\
\mathcal{E} =  \{(v,v')\in \mathcal{V}\times \mathcal{V} :
|x-x'|=1)\},
\end{array}
\]
where $|v-v'|=\sum_{i=1}^d|v(i)-v'(i)|$ is the Hamming distance in
$\mathcal{V}$.

We will consider symmetric trap models in $\H$, namely
continuous time, space inhomogeneous, simple random walks in $\H$, whose 
transition probabilities (from each site of $\H$ to any of its $d$ nearest neighbors) 
are uniform. Let $\g^{d}=\{\g^{d}_v,\,v\in\V\}$ denote the set of mean waiting
times characterizing the model.

We will map $\V$ onto the set $\D:=\{1,\ldots,2^d\}$ by enumerating $\V$
in decreasing order of $\g^{d}$ (with an arbitrary tie breaking
rule), and then consider $X^{d}$, the mapped process. Let
\begin{equation}\label{eq:tgd}
\tgd=\{\tg^{d}_x,\,x\in\D\}
\end{equation}
denote the enumeration in decreasing order of $\g^{d}$, and view
it as a finite measure in $\N=\{1,2,\ldots\}$, the positive
natural numbers.

We next consider a class of processes which turns out to contain limits of trap models
in $\H$ as $d\to\infty$, as we will see below. Let  ${\cal
N}:=\{(N^{(x)}_t)_{t\geq0},\,x\in\N\}$ be i.i.d.~Poisson processes
of rate 1, with $\s^{(x)}_j$ the $j$-th event time of $N^{(x)}$, $x\in\N$, $j\geq1$,
and let ${\cal T}=\{T_0;\,T^{(x)}_i\,i\geq1,x\in\N\}$ be
i.i.d.~exponential random variables of rate 1. ${\cal N}$ and
${\cal T}$ are assumed independent. Consider now a finite measure
$\gamma$ supported on $\N$, and for $y\in\bN=\N\cup\{\infty\}$ let
\begin{equation}
  \label{eq:G}
  \Gamma(t)=\Gamma^{y}(t)
  =\g_y\,T_0+\sum_{x=1}^\infty\g_x\sum_{i=1}^{N^{(x)}_t}T^{(x)}_i,
\end{equation}
where, by convention, $\sum_{i=1}^{0}T^{(x)}_i=0$ for every $x$,
and $\gamma_\infty=0$. We define the process $Y$ on $\bN$ starting
at $y\in\bN$ as follows. For $t\geq0$
\begin{equation}
  \label{eq:X}
  Y_t=\begin{cases}\mbox{}\,y,&\mbox{ if } 0\leq t<\g(y)\,T_0,\\
          \mbox{}\,x,&\mbox{ if } \Gamma(\s^{(x)}_j-)\leq
          t<\Gamma(\s^{(x)}_j)\mbox{ for some }1\leq j<\infty,\\
          \infty,&\mbox{ otherwise}.\end{cases}
\end{equation}
This process, which we here call the K process with parameter
$\gamma$, was introduced and studied in~\cite{kn:FM}, where it was shown to
arise as limits of trap models in the complete graph with $n$
vertices as $n\to\infty$ (see Lemma 3.11 in~\cite{kn:FM}).
In the next section, we derive
a similar result for the hypercube. See Theorem~\ref{teo:main}.
This is our main technical result. Then, in the following section,
we apply that result to get the scaling limits of the REM-like
trap model and the Random Hopping Times dynamics for the REM in
ergodic time scales as K processes. See Section~\ref{sec:brht}.

\section{Convergence to the K process}
\setcounter{equation}{0}
\label{sec:conv}

\begin{theo}\label{teo:main}
Suppose that, as $d\to\infty$, $\tgd$ converges weakly to a finite
measure $\tg$ supported on $\N$, and that $X^d_0$ converges weakly
to a probability measure $\mu$ on $\bN$. Then, $X^{d}$ converges
weakly in Skorohod space as $d\to\infty$ to a K process with
parameter $\tg$ and initial measure $\mu$.
\end{theo}

This result extends the analysis performed in~\cite{kn:FM} for the trap
model in the complete graph, with a similar approach (see Lemma
3.11 in~\cite{kn:FM} and its proof). The extra difficulty here
comes from the fact that the transition probabilities in the
hypercube are not uniform in the state space, as is the case in
the complete graph. However, all that is indeed needed is an
approximate uniform entrance law in finite sets of states.  This
result, a key tool used several times below, is available
from~\cite{kn:AG}. We state it next, in a form suitable to our
purposes, but first some notation. Let $\X$ denote the embedded
chain of $X^d$ and for a given fixed finite subset $\J$ of $\N$,
let $\T_{\J}$ denote the entrance time of $\X$ in $\J$, namely,
\begin{equation}\label{eq:tj}
\T_{\J}=\inf\{n\geq0:\X_n\in\J\}.
\end{equation}
\begin{prop}[Corollary 1.5~\cite{kn:AG}]\label{prop:ag}
\begin{equation}\label{eq:un}
\lim_{d\to\infty}\max_{x\notin\J,y\in\J}\left|\P(\X_{\T_\J}=y|\X_0=x)-\frac1{|\J|}\right|=0.
\end{equation}
\end{prop}
Here $|\cdot|$ denotes cardinality.
\begin{rmk}\label{rmk:ag}
Corollary 1.5 of~\cite{kn:AG} is actually more precise and
stronger than the above statement, with error of approximation
estimates, and holding for $\J$ depending on $d$ in a certain
manner as well.
\end{rmk}

\begin{rmk}\label{rmk:hc}
(\ref{eq:un}) is the only fact about the hypercube used in the proof of Theorem~\ref{teo:main}.
This result would thus hold as well for other graphs with the same property. The hypercube has nevertheless 
been singled out in analyses of dynamics of mean field spin glasses (see above mentioned
references), and that is a reason for us to do the same here.
\end{rmk}

\subsection{Proof of Theorem~\ref{teo:main}}
\label{ssec:proof}

The strategy is to approximate $X^d$ for $d$ large by a trap model
in the complete graph with vertex set $\M=\{1,\ldots,M\}$ and mean
waiting times $\{\tg_1,\ldots,\tg_M\}$ for $M\leq d$ large.
Let $Y^M$ denote the latter process, and let us put
$Y^M_0=Y_0\,1\{Y_0\in\M\}+W\,1\{Y_0\notin\M\}$,
with $W$ an independent uniform random variable in $\M$. To accomplish the
approximation, we will resort to an intermediate process, which we
next describe. We start by considering $\hat X^d$, the trap model
on $\D$ obtained from $X^d$ by replacing its set of
mean waiting times (see~(\ref{eq:tgd}) above) by
$\{\tg_x,\,x\in\D\}$. The intermediate process we will
consider is then $\hat X^d$ {\em restricted} to $\M$, denoted
$\hat X^{d,M}$: this is the Markov process obtained from $\hat X^d$ by
observing it only when it is in $\M$ (with time stopping for
$\hat X^{d,M}$ when $\hat X^d$ is outside $\M$).

The approximations will be strong ones: we will couple $X^d$ to
$\hat X^{d,M}$ and $\hat X^{d,M}$ to $Y^M$, in the spirit of  Theorem 5.2
in~\cite{kn:FM}, where the approximation of $Y$ by $Y^M$, needed
here as the last step of the argument, was established. In particular, we also couple $X_0^d$ to
$Y_0$ so that the former converges almost surely to the latter as $d\to\infty$.

\subsubsection{Coupling of $\hat X^{d,M}$ and $Y^M$}
\label{sssec:coup1}

We first look at the embedded chains of $\hat X^{d,M}$ and $Y^M$. Let
$(p^{d,M}_{ij})_{i,j\in\M}$ be the transition probabilities of the
former chain, and let $\hat p=\min_{i,j\in\M}p^{d,M}_{ij}$. We
leave it to the reader to check that there is a coupling between
both chains which agrees at each step with probability at least
$M\hat p$. We resort to such a coupling. Proposition~\ref{prop:ag}
implies that
\begin{equation}\label{eq:mps}
    M\hat p\to1
\end{equation}
as $d\to\infty$ for every $M$ fixed.

Since $\hat X^{d,M}$ and $Y^M$
have the same mean waiting times at each
site, we can couple them in such a way that they have the same
waiting times at successive visits to each site.
One may also find a coupling of $\hat X^{d,M}_0$ and $Y^M_0$ such that
\begin{equation}\label{eq:ic}
\P(\hat X^{d,M}_0\ne Y^M_0)\to0
\end{equation}
as $d\to\infty$. Resorting also to that
coupling, we get the following result.

\begin{lm}\label{prop:coup1}
For every $T$ and $M$ fixed, we have
\begin{equation}\label{eq:c1}
    \P(\hat X^{d,M}_t=Y^M_t,\,t\in[0,T])\to1
\end{equation}
as $d\to\infty$.
\end{lm}

\noindent{\bf Proof}

Let $N_T$
denote the number of jumps of $Y^M$ in the time interval $[0,T]$,
and $0=t_0,t_1,\ldots,t_{N_T}$ the respective jump times. We conclude
from the above discussion that
\begin{eqnarray}\nonumber
\P(\hat X^{d,M}_t= Y^{M}_t, t\in[0,T]| N_T=k) &=& \P(\hat X^{d,M}_{t_i}=
Y^{M}_{t_i}, i=0,1,\ldots,k| N_T=k)\\
\label{eq:c2}&\geq&\P(\hat X^{d,M}_0=Y^M_0)\,(M\hat p)^k,
\end{eqnarray}
and the result follows from~(\ref{eq:mps}),~(\ref{eq:ic}) and
dominated convergence (since $M\hat p$ is bounded above by 1 for
all $d$ and $M$). $\square$

\subsubsection{Coupling of $X^d$ and $\hat X^{d,M}$}
\label{sssec:coup2}

We couple $X^d$ and $\hat X^d$ (the latter process was defined at
the beginning of the section) in the following way. Note that $\X$
is their common embedded chain. We then make the successive
sojourn times of $X^d$ and $\hat X^d$, starting from the first
ones, be given by
$\tg^d_{\X_1}T^{\X_1}_1,\tg^d_{\X_2}T^{\X_2}_2,\ldots$ and
$\tg_{\X_1}T^{\X_1}_1,\tg_{\X_2}T^{\X_2}_2,\ldots$, respectively,
where the common $T^x_1,T^x_2,\ldots$, $x\in\N$, are i.i.d.~rate 1
exponential random variables.

With a view towards approximating $X^d$ and $\hat X^{d,M}$ strongly in Skorohod space, we
introduce a {\em time distortion function} useful for that (see~(\ref{eq:la}) below). For $K$ a fixed positive integer, let
$\K$ denote the set $\{1,\ldots,K\}$, and consider the successive entrance and exit times of
$X^d$ and $\hat X^{d,M}$ in and out of $\K$ defined as follows. Let $\tau_0=\tau^*_0=\xi_0=\xi^*_0=0$ and for $i\geq1$,
let
\begin{eqnarray}\label{eq:tau}
\tau_i=\inf\{t\geq\tau^*_{i-1}:\,\hat X^{d,M}_t\in\K\},\quad
\tau^*_i=\inf\{t\geq\tau_{i}:\,\hat X^{d,M}_t\notin\K\},
\end{eqnarray}
and similarly define $\xi_i$ and $\xi^*_i$, $i\geq1$ with $X^{d}$ replacing $\hat X^{d,M}$.
See Figure~\ref{fig:1} below.

\begin{figure}[!htb]
        \centering
                \includegraphics[scale=0.635]{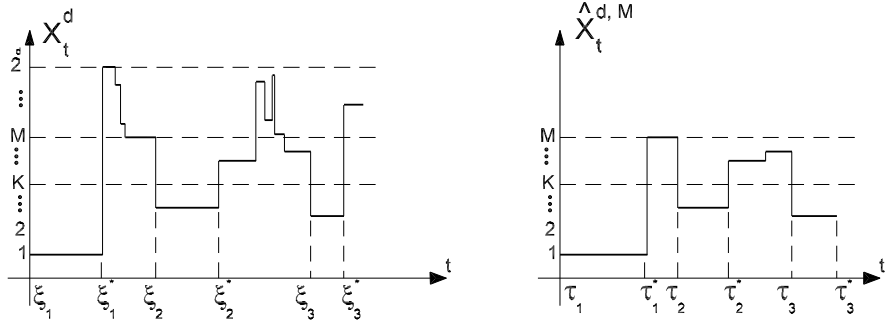}
        \caption{Trajectories of $X^d$ and $\hat X^{d,M}$.}
        \label{fig:1}
\end{figure}

Now for $T>0$ fixed, let $N=\min\{i\geq1:\,\tau_i>T\}$ and define
\begin{equation}\label{eq:la}
\tl_t=\begin{cases}
\xi_j +\frac{\xi^*_j-\xi_j}{\tau^*_j-\tau_j}\,(t-\tau_j), & \textrm{if} \ \tau_j<t\leq \tau_j^{*}\ \textrm{for some $0\leq j<N$},\\
\xi^{*}_j+\frac{\xi_{j+1}-\xi_{j}^{*}}{\tau_{j+1}-\tau_j^{*}}\,(t-\tau^{*}_j),
& \textrm{if}  \ \tau^{*}_j<t\leq \tau_{j+1}\ \textrm{for some $0\leq j<N$},\\
\xi_{N} +(t-\tau_{N}),& \textrm{if} \ t\geq \tau_N.
\end{cases}
\end{equation}
Here and below, we interpret $0/0$ as $1$. See Figure~\ref{fig:2} below.

\begin{figure}[!htb]
        \centering
                \includegraphics[scale=0.635]{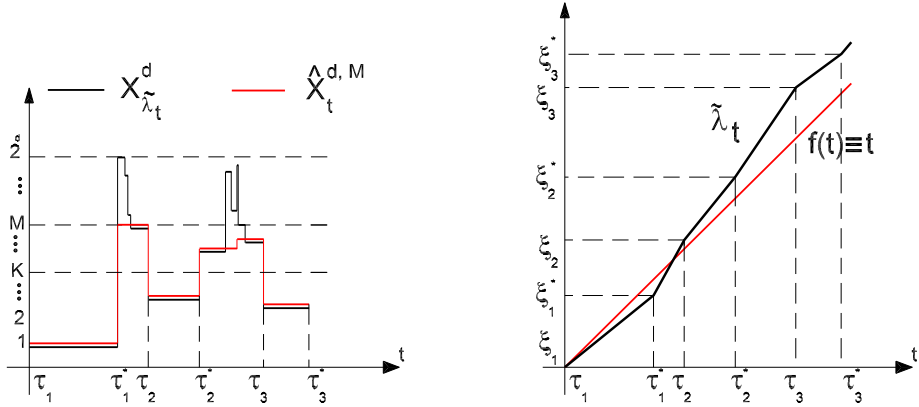}
        \caption{Superimposed trajectories of $X^d$ with time distorted by $\tilde\lambda$ and $\hat X^{d,M}$ (left), and superposition of the graphs
                 of $\tl$ and the identitity (right).}
        \label{fig:2}
\end{figure}

\begin{rmk}\label{rmk:coup2}
With the above definition of $\tl$, we first note that $\hat X^{d,M}_t=X^d_{\tilde\l_t}$ whenever any of both processes is visiting $\K$ before time $\tau_N$.
\end{rmk}

As part of the norm in Skorohod space, we consider the class $\L$
of nondecreasing Lipschitz functions mapped from $[0,\infty)$ onto
$[0,\infty)$, and the following function on $\L$.
\begin{equation}\label{eq:g}
    \phi(\l)=\sup_{0\leq s<t}\left|\log\frac{\l_t-\l_s}{t-s}\right|
\end{equation}
We have from~(\ref{eq:la}) above that
\begin{equation}\label{eq:g1}
    \phi(\tl)\leq\max_{1\leq j\leq N}\left|\log\frac{\xi_{j}-\xi_{j-1}^{*}}{\tau_{j}-\tau_{j-1}^{*}}\right|
    \vee
    \max_{1\leq j\leq N}\left|\log\frac{\xi^*_j-\xi_j}{\tau^*_j-\tau_j}\right|.
\end{equation}

Below, we will consider the events $A_j$, $j\geq0$, as follows.
\begin{eqnarray}\label{eq:B}
    A_0\!\!\!\!\=\!\!\!\!\{X^d_0\in\K\}\cup\{\mbox{there exists $0\leq t<\xi_{1}$ such that }
      X^d_t\in\M\setminus\K\},\\
    A_j\!\!\!\!\=\!\!\!\!\{\mbox{there exists $\xi^*_j\leq t<\xi_{j+1}$ such that }
      X^d_t\in\M\setminus\K\},\,j\geq1.
\end{eqnarray}
It follows from Proposition~\ref{prop:ag} that for $j\geq0$
\begin{equation}\label{eq:b1}
\lim_{d\to\infty}\inf_{x\notin\M}\P(A_j|X^d_{\xi_{j}^*}=x)=1-K/M.
\end{equation}
Notice that the probability on the left hand side of~(\ref{eq:b1})
does not depend on $j$; we thus get that
\begin{equation}\label{eq:b2}
\lim_{M\to\infty}\lim_{d\to\infty}\P(A_j)=1 \mbox{ uniformly on }
j\geq0.
\end{equation}

\subsubsection{Conclusion of proof of Theorem~\ref{teo:main}}
\label{sssec:proof}

Let $D_{\bN}([0,\infty))$ be the (Skorohod) space of c\'adl\'ag functions of
$[0,\infty)$ to $\bN$ with metric
\begin{equation}\label{eq:met1}
\rho(f,g) := \inf_{\l\in\Lambda} \left[\phi(\l) \vee \int_0^\infty
e^{-u} \rho(f,g,\l,u) du \right],
\end{equation}
where
\begin{equation}\label{eq:met2}
\rho(f,g,\l,u) := \sup_{t\geq 0}\left|[f(t\wedge
u)]^{-1}-[g(\l(t)\wedge u)]^{-1}\right|;
\end{equation}
see Section 3.5 in~\cite{kn:EK}; $\L$ and $\phi$ were defined in the paragraph
of~(\ref{eq:g}) above, and $\infty^{-1}=0$.

It follows from Lemma 3.11 in~\cite{kn:FM} that $Y^M$ converges
weakly to $Y$ in Skorohod space as $M\to\infty$. In order to prove
Theorem~\ref{teo:main}, it thus suffices to show the following
result.

\begin{lm}\label{prop:main}
With the above construction of $X^d$ and $Y^M$, we have that for every
$\e>0$
  \begin{equation}
    \label{eq:prob}
   \lim_{M\to\infty}\limsup_{d\to\infty}\P(\rho(X^d,Y^M)>\e)=0.
  \end{equation}
\end{lm}

\noindent{\bf Proof}

Given $\e>0$, let $T_\e=-\log(\e/2)$. Then choosing $\l$ to be the
identity, noticing that $\rho$ in~(\ref{eq:met2}) is bounded above by
$1$, and using Lemma~\ref{prop:coup1}, we find that for
every $M>0$
  \begin{equation}
    \label{eq:prob1}
   \P(\rho(\hat X^{d,M},Y^M)>\e/2)\leq\P(\hat X^{d,M}_t\ne Y^M_t\mbox{ for some
   }t\in[0,T_\e])\to0
  \end{equation}
as $d\to\infty$. So, to establish Lemma~\ref{prop:main}, it
suffices to prove Lemma~\ref{lm:main} below.
$\square$

\begin{lm}\label{lm:main}
With above construction of $X^d$ and $\hat X^{d,M}$, we have that for
every $\e>0$
  \begin{equation}
    \label{eq:prob2}
   \lim_{M\to\infty}\limsup_{d\to\infty}\P(\rho(\hat X^{d,M},X^d)>\e)=0.
  \end{equation}
\end{lm}

\noindent{\bf Proof}

Let $T=T_\e=-\log\e$, choose $K$ such that $|x^{-1}-y^{-1}|\leq\e$
for every $x,y\in\bN\setminus\K$, and consider $\tl$ as
in~(\ref{eq:la}) with such $T$ and $K$. Then, by
Remark~\ref{rmk:coup2} and~(\ref{eq:g}), we see that it suffices
to show that for every $\e>0$
  \begin{equation}
    \label{eq:prob4}
   \lim_{M\to\infty}\limsup_{d\to\infty}
   \P\left(\max_{1\leq j\leq N}\left|\log\frac{\xi^*_j-\xi_j}{\tau^*_j-\tau_j}\right|>\e\right)=0,
  \end{equation}
and
  \begin{equation}
    \label{eq:prob3}
   \lim_{M\to\infty}\limsup_{d\to\infty}
   \P\left(\max_{1\leq j\leq N}\left|\log\frac{\xi_{j}-\xi_{j-1}^{*}}{\tau_{j}-\tau_{j-1}^{*}}\right|>\e\right)=0.
  \end{equation}

{\bf Proof of~(\ref{eq:prob4})}

One readily checks that
  \begin{equation}
    \label{eq:p41}
   \max_{1\leq j\leq N}\left|\log\frac{\xi^*_j-\xi_j}{\tau^*_j-\tau_j}\right|
   \leq
   \max_{x\in\K}\left|\log\frac{\tgd_{x}}{\tg_{x}}\right|,
 \end{equation}
and~(\ref{eq:prob4}) follows immediately from the assumption that
$\tgd\to\tg$ as $d\to\infty$.

\medskip

{\bf Proof of~(\ref{eq:prob3})}

Let $\tg^{d,K}=\{\tg^{d,K}_x:=\tg^d_x\wedge\tg^d_K,x\in\D\}$, and consider the
trap model $X^{d,K}$ with mean waiting times $\tg^{d,K}$ coupled
to $X^{d}$ so that both processes have the same embedded chain
$\X$ and the respective sojourn times are given by
$\tg^{d,K}_{\X_1}T^{\X_1}_1,\tg^{d,K}_{\X_2}T^{\X_2}_2,\ldots$ and
$\tg^d_{\X_1}T^{\X_1}_1,\tg^d_{\X_2}T^{\X_2}_2,\ldots$

Let now $\tilde X^{d,K}$ denote the process $X^{d,K}$ restricted
to $\K$ (analogously as $\hat X^{d,M}$), with $\tilde\X$ its embedded
chain. Let $N^K$ denote the number of jumps of $\tilde X^{d,K}$ up
to time $T$. Notice that $N^K$ is a Poisson process with rate
$1/\tgd_K$ independent of $\X$ and of the history of $X^d$ in the
time intervals $[\xi^*_j,\xi_{j+1})$, $j\geq0$. Thus, the
probability on the left hand side of~(\ref{eq:prob3}) is bounded
above by
  \begin{equation}
    \label{eq:p31}
   \P\left(\max_{1\leq j\leq
   N^K}\left|\log\frac{U_j}{V_j}\right|>\e\right)\leq\sum_{n=1}^\infty\sum_{j=1}^n
   \P\left(\left|\log\frac{U_j}{V_j}\right|>\e\right)\P(N^K=n),
  \end{equation}
where $U_j:=\xi_{j}-\xi_{j-1}^{*}$ and
$V_j:=\tau_{j}-\tau_{j-1}^{*}$. We now estimate the first
probability on the right hand side of~(\ref{eq:p31}). We first
note that from~(\ref{eq:b2}), and the fact that $\E(N^K)$ is
finite and independent of $d,M$, we may insert $A_j$ in that
probability. We next write $U_j=W_j+R_j$, where $R_j$ is the time
spent by $X^d$ in $\D\setminus\M$ during the time interval
$[\xi^*_{j-1},\xi_j)$. From the elementary inequality
$|\log(x+y)|\leq|\log x|+y$, valid for all $x,y>0$, we get that
  \begin{equation}
    \label{eq:p32}
   \P\left(\left|\log\frac{U_j}{V_j}\right|>\e,\,A_j\right)\leq
   \P\left(\left|\log\frac{W_j}{V_j}\right|>\e/2,\,A_j\right)+
   \P\left(R_j>\e V_j/2\right).
  \end{equation}
Arguing as in~(\ref{eq:p41}) above, we find that the first event
in the first probability on the right hand side of~(\ref{eq:p32})
is empty as soon as
$\max_{x\in\M}\left|\log\frac{\tgd_{x}}{\tg_{x}}\right|\leq\e/2$,
thus from $\tgd\to\tg$ as $d\to\infty$ we only need to consider the
second probability on the right hand side of~(\ref{eq:p32}). One
readily checks that it is bounded above by
  \begin{equation}
    \label{eq:p33}
   \max_{x\notin\K}\P\!\left(R_j>\e V_j/2\left|X^d_{\xi^*_{j-1}}=x\right.\right)
  \end{equation}
for all $j\geq0$, and that the above expression does not depend on
$j$. It is enough then to show that for any $\e>0$
  \begin{equation}
    \label{eq:p34}
   \P\!\left(R_1>\e V_1|X^d_0=x\right)=:
   \P_x(R_1>\e V_1)\to0
  \end{equation}
as $d\to\infty$ and then $M\to\infty$, uniformly in $x>K$. This is
readily seen to follow from the facts that
  \begin{equation}
    \label{eq:p35}
  \lim_{M\to\infty}\limsup_{d\to\infty}\max_{x\notin\K}\P_x(R_1>\e)=0
  \end{equation}
for any $\e>0$, and that, given $\d>0$, there exists $\e>0$ such that
  \begin{equation}
    \label{eq:p36}
  \limsup_{M\to\infty}\limsup_{d\to\infty}\max_{x\notin\K}\P_x(V_1\leq\e)\leq\d.
  \end{equation}

{\bf Proof of~(\ref{eq:p35})}

Let $x\notin\K$ be arbitrary. We will estimate
\begin{equation}\label{eq:p38}
\E_x(R_1):=\E(R_1|X^d_0=x)=\sum_{y=M+1}^d\tgd_y\,\E_x(\ll(y)),
\end{equation}
where $\ll(y)$ is the number of visits of $\X$ to $y$ from time $0$ till its first entrance in $\K$.

Let $\K_y=\K\cup\{y\}$  and consider the discrete time Markov
process $\bar\X$ such that $\bar\X_0=x$ and otherwise $\bar\X$ is
the restriction of $\X$ to $\K_y=\K\cup\{y\}$, and let
$\bar\ll(y)$ denote the number of visits of $\bar\X$ to $y$ from
time $0$ till its first entrance in $\K$.  Clearly,
\begin{equation}\label{eq:p40}
\ll(y)=\bar\ll(y).
\end{equation}
Now let $\X^*$ denote the Markov chain on $\K_y\cup\{x\}$ with the
following set of transition probabilities. Let
$p^{1}=(p^{1}_{wz},\,w,z\in\K_y\cup\{x\})$, and
$p^{2}=(p^{2}_{wz},\,w,z\in\K_y\cup\{x\})$ denote the sets of
transition probabilities of $\bar \X$ and $\X^*$, respectively. We
make
\begin{equation}\label{eq:p39}
p^{2}_{xy}=p^{2}_{yy}=p^*:=\max\{p^{1}_{wz};\,w=x,y;\,z\in\K_y\},
\end{equation}
and the remaining $p^{2}_{wz}$ can be assigned arbitrarily with the
only obvious condition that $p^{2}$ is a set of transition
probabilities on $\K_y$. Let $\ll^*(y)$ denote the number of visits of
$\X^*$ to $y$ from time $0$ till its first entrance in $\K$.
One readily checks that $\ll^*(y)$ is a Geometric random variable
with parameter $1-p^*$ and that it stochastically dominates
$\bar\ll(y)$. From this and~(\ref{eq:p40}), we conclude that
\begin{equation}\label{eq:p42}
\E_x(\ll(y))\leq\frac{p^*}{1-p^*}
\end{equation}
uniformly in $x\notin\K$.  Proposition~\ref{prop:ag} then implies
that
\begin{equation}\label{eq:p43}
\limsup_{d\to\infty}\max_{x\notin\K}\E_x(\ll(y))\leq\frac{1}{K}.
\end{equation}
It follows readily from this,~(\ref{eq:p38}) and the assumption
that $\tgd\to\tg$ as $d\to\infty$ that
\begin{equation}\label{eq:p44}
\limsup_{d\to\infty}\max_{x\notin\K}\E_x(R_1)\leq\frac{1}{K}\sum_{y=M+1}^\infty\tg_y,
\end{equation}
and~(\ref{eq:p35}) follows (using Markov's inequality), since $\tg$
is a finite measure on $\N$.

\medskip

{\bf Proof of~(\ref{eq:p36})}

Let us fix $x_0\notin\K$. Consider the Markov process $Z=(Z_t)_{t\geq0}$ on $\M$ such that $Z_0=x_0$, for every $x\in\M$ the waiting time at $x$ before jumping is exponential with mean $\tg_x$, and the transition probability to $y\in\M$ equals $\hat p$, if $y\notin\K$, and $\frac{1-(M-K)\hat p}{K}$,
if $y\in\K$, where $\hat p$ was defined in the paragraph of~(\ref{eq:mps}) above.

One readily checks that, given $X^d_0=x_0$, $V_1$ stochastically dominates $S$, the time $Z$ spends in $\M\setminus\K$ before hitting $\K$ for the first time. Since $\tg_x$ is decreasing in $x$, by the construction of $Z$, we have that, for every $L\in\{K+1,\ldots,M\}$, $S$ dominates stochastically the random variable $\tg_LT1_C$, where $C$ is the event that $Z$ visits $\{K+1,\ldots,L\}$ before hitting $\K$ for the first time, $T$ is an exponential random variable of rate 1, and $1_{\cdot}$ is the usual indicator function.

Let now $\zz$ denote the embedded chain for $Z$, and $\T_L=\inf\{n\geq1:\zz_n\in\{1,\ldots,L\}\}$. Then
\begin{eqnarray}\nonumber
&&\P(C|\zz_0=x_0)\\\nonumber\=\P(\zz_{\T_L}\in\{K+1,\ldots,L\}|\zz_{\T_L}\leq L,\zz_0=x_0)\\\nonumber
&=&\frac{(L-K)\hat p}{K\left[\frac{1-(M-K)\hat p}{K} \right] + (L-K)\hat p} = \frac{(L-K)\hat p}{1-(M-L)\hat p}\\
&=&\frac{(L-K)M\hat p}{(1-M\hat p)M+LM\hat p}\to\frac{(L-K)}{L}=1-\frac KL
\end{eqnarray}
as $d\to\infty$ uniformly in $x_0$; see~(\ref{eq:mps}).

We conclude that
  \begin{equation}
    \label{eq:p37}
\limsup_{M\to\infty}\limsup_{d\to\infty}\max_{x\notin\K}\P_x(V_1\leq\e)\leq\P(\tg_LT\leq\e)+\frac KL
  \end{equation}
for every $K<L$. Thus, given $K$ and $\d>0$, we first choose $L$ such that $\frac KL\leq\d/2$, and then $\e>0$
such that $\P(\tg_LT\leq\e)\leq\d/2$, and we are done.
$\square$

\section{The REM-like trap model and the Random hopping times dynamics for the REM}
\setcounter{equation}{0}
\label{sec:brht}

In this section we apply Theorem~\ref{teo:main} to obtain the scaling limits of two
disordered trap models in the hypercube, namely trap models in the hypercube whose
mean waiting times are random variables.

\subsection{The REM-like trap model}
\setcounter{equation}{0}
\label{ssec:btm}

Let $\tau^{d}:=\{\tau^{d}_v,\,v\in\V\}$, an i.i.d.~family of random variables in the domain
of attraction of an $\a$-stable law with $0<\a<1$, be the mean waiting times of a trap
model in $\H$. Let us then consider as before $\ttau^{d}:=\{\ttau^{d}_x,\,x\in\D\}$,
the decreasing order statistics of $\tau$ (with an arbitrary tie breaking
rule), and let $Y^{d}$ be the mapped process on $\D$.

Now let $c_d$ be a scaling factor such that $\htau^{d}:=c_d\,\ttau^{d}$ converges to
the increments in $[0,1]$ of an $\a$-stable subordinator. Namely,
\begin{equation}
\label{eq:cd}
c_d=\left(\inf\{t\geq0:\P(\tau_0>t)\leq 2^{-d}\}\right)^{-1}.
\end{equation}

Let us now consider $Y^{d}$ speeded up by $c_d$, namely $\hy_t=Y^{d}(t/c_d)$, $t\geq0$.
Notice that $\hy$ is a trap model on $\H$ with mean waiting times given by $\htau^{d}$.

Let $\hat\gamma=\{\hat\g_i,\,i\in\N\}$ denote the increments in $[0,1]$ of an $\a$-stable subordinator
in decreasing order.

\begin{cor}\label{cor:btm}
Suppose that $\htau^{d}_0$ converges weakly to a probability
measure $\mu$ on $\bN$. Then
  \begin{equation}
    \label{eq:btm}
    (\hy,\htau^{d})\Rightarrow(Y,\hat\g),
  \end{equation}
where $Y$ is a K process with parameter $\hat\g$ and initial measure
$\mu$, and $\Rightarrow$ means weak convergence in the product of
Skorohod norm and weak convergence norm in the space of finite
measures in $\N$.
\end{cor}

\begin{rmk}\label{rmk:fm}
In~\cite{kn:FM}, a similar result was proved for the REM-like trap model in the complete graph.
See Theorem 5.2 in that reference.
\end{rmk}

\noindent{\bf Proof of Corollary~\ref{cor:btm}}

We can suppose that we are in a probability space where $\htau^{d}\to\hat\g$ almost surely
(see proof of Theorem 5.2 in~\cite{kn:FM} for an explicit argument). We can then invoke
Theorem~\ref{teo:main} to get that $\hy\Rightarrow Y$, and the full result follows. $\square$

\subsection{Random Hopping Times dynamics for the REM}
\setcounter{equation}{0}
\label{ssec:rht}

This is a dynamics whose equilibrium is the Random Energy Model. Let
$H^{d}:=\{H^{d}_v,\,v\in\V\}$ be an i.i.d.~family of standard normal random variables, and make
$\tau^{d}:=\{\tau^{d}_v,\,v\in\V\}$, where $\tau^{d}_v=e^{\b{\sqrt d}H^{d}_v}$. Defining now
$\ttau^{d}$, $\htau^{d}$, $Y^{d}$ and $\hy$ as above, with
\begin{equation}
\label{eq:cd1}
c_d=e^{-\frac{2\log 2}\a d+\frac1{2\a}\log d},
\end{equation}
with $\a=\sqrt{2\log 2}/\b$, we have that, if $\a<1$, then
Corollary~\ref{cor:btm} holds in this case as well, with $\hat\g$ and
$\mu$ as before.

The proof starts from the known result that in this case $\htau^{d}\Rightarrow\hat\g$ (see Remark to Theorem 2
in~\cite{kn:GPM}). Again, as in the proof of Corollary~\ref{cor:btm} above, we can go to a probability space
where the latter convergence is almost sure, and close the argument in the same way.

\begin{rmk}\label{rmk:erg}
The time scale $t\to t/c_d$ adopted in the above models is the {\em ergodic} time scale mentioned in~\cite{kn:BD}.
Under shorter scalings (i.e., $t\to t/c'_d$, with $c'_d>\!\!>c_d$) the model exhibits {\em aging} (when starting
from the uniform distribution),
and under longer ones ($c'_d<\!\!<c_d$), the model reaches equilibrium. More precisely, under shorter scalings,
we have that as $d\to\infty$
\begin{equation}\label{eq:age1}
 \P_{\mu^d}\left(Y^{d}(t/c'_d)=Y^{d}((t+s)/c'_d)\right)\to{\cal R}(s/t),
\end{equation} 
where $\mu^d$ is the initial uniform distribution on
$\D$, and ${\cal R}$ is a nontrivial function such
that ${\cal R}(0)=1$ and $\lim_{x\to\infty}{\cal R}(x)=0$. Indeed,
for the models of this section (as well as in many other instances
in the references), ${\cal R}$ is the arcsine law: 
\begin{equation}\label{eq:arcsin}
{\cal R}(x)=\frac{\sin(\pi\a)}{\pi}\int_{\frac{x}{1+x}}^1 s^{-\a}(1-s)^{\a-1}\,ds.
\end{equation} 
See~\cite{kn:AC2}. Under longer scalings, it can be shown that $Y^{d}(t/c'_d)\Rightarrow\bar\g$ as $d\to\infty$ for every $t>0$,
where $\bar\g$ is $\hat\g$ normalized to be a probability measure. It is the limiting equilibrium measure, or more precisely,
the equilibrium measure of $Y$.
\end{rmk}

\begin{rmk}\label{rmk:age}
It can be shown that $Y$ exhibits aging at a vanishing time scale (when starting from $\infty$), i.e.
\begin{equation}\label{eq:age2}
\P_\infty\left(Y(\epsilon t)=Y(\epsilon (t+s))\right)\to{\cal R}(s/t)
\end{equation}
as $\epsilon\to0$. See Theorem 5.11 in~\cite{kn:FM}. This is in agreement with~(\ref{eq:age1}).
\end{rmk}

\section*{Acknowledgements}
This work was mostly done as part of the master's project of the second author, at IME-USP and
financed by FAPESP. We thank Claudio Landim for pointing out a mistake in a previous version of this paper.

\end{document}